\documentclass[draft]{amsart}

\newtheorem{thm}{Theorem}
\newtheorem{prop}[thm]{Proposition}

\theoremstyle{remark}
\newtheorem{rmk}[thm]{Remark}

\numberwithin{thm}{section}
\numberwithin{equation}{section}

\newcommand{\nv}{\mathbf{N}}
\newcommand{\rv}{\mathbf{R}}

\begin{document}

\title[The Briot-Bouquet systems and the center families]{The Briot-Bouquet systems and the center families for holomorphic dynamical systems}

\author[F. Rong]{Feng Rong}

\address{Department of Mathematics, Shanghai Jiao Tong University, 800 Dong Chuan Road, Shanghai, 200240, P.R. China}

\email{frong@sjtu.edu.cn}



\thanks{The author is partially supported by the National Natural Science Foundation of China (Grant No. 11001172), the Specialized Research Fund for the Doctoral Program of Higher Education of China (Grant No. 20100073120067) and the Scientific Research Starting Foundation for Returned Overseas Chinese Scholars.}

\begin{abstract}
We give a complete solution to the existence of isochronous center families for holomorphic dynamical systems. The study of center families for $n$-dimensional holomorphic dynamical systems naturally leads to the study of ($n-1$)-dimensional Briot-Bouquet systems in the phase space. We first give a detailed study of the Briot-Bouquet systems. Then we show the existence of isochronous center families in the neighborhood of the equilibrium point of three-dimensional systems based on the two-dimensional Briot-Bouquet theory. The same approach works in arbitrary dimensions.
\end{abstract}

\maketitle

\section{Introduction}

We consider the holomorphic dynamical system
\begin{equation}\label{E:I}
\dot{z_i}=F_i(z_1,\cdots,z_n),\ \ \ \ \ \ (i=1,\cdots,n),
\end{equation}
where $F_i$'s are holomorphic functions with $F_i(0,\cdots,0)=0$, i.e. the origin is an equilibrium point of (\ref{E:I}). In particular, we are interested in the existence of isochronous center families in the neighborhood of the origin.

Let $\Lambda (z_1,\cdots,z_n)^T$ denote the linear part of $(F_1,\cdots,F_n)$. The eigenvalues of $\Lambda$ are determined up to a real scaling constant (due to the scaling of the time variable). By the well-known stable/unstable manifold theorem, $\Lambda$ has at least one purely imaginary eigenvalue if a center family exists in the neighborhood of the origin. It is also well-known that if all the eigenvalues of $\Lambda$ are purely imaginary and equal and $\Lambda$ is diagonalizable, then the origin is an isochronous center (due to a holomorphic linearization in the neighborhood of the origin, cf. \cite{A}).

In \cite{N}, Needham showed the existence of an isochronous center family in the neighborhood of the origin for (\ref{E:I}) assuming that $\Lambda$ has only one purely imaginary eigenvalue. In \cite{NM}, Needham and McAllister focused on two-dimensional systems and proved similar results assuming that both of the eigenvalues of $\Lambda$ are purely imaginary and $\Lambda$ is diagonalizable. In both papers, the approach is to study (\ref{E:I}) in the phase space and apply the Briot-Bouquet singular point theory. (Some results on isochronous center families have also been obtained by Zhang (\cite{Z1,Z2}), using other methods.)

In this paper, we use a similar approach to study (\ref{E:I}) in arbitrary dimensions. For this purpose, we need a detailed knowledge of the Briot-Bouquet systems in general. We give such a detailed study in Section \ref{S:BB}. As an illustration, we give a complete description of isochronous center families in the neighborhood of the equilibrium point of three-dimensional systems in Section \ref{S:CF}. One can obtain similar description for arbitrary dimensional systems using the same method.

Note that in our study we do not assume $\Lambda$ to be diagonalizable (cf. Remark \ref{R:C2}). We also note that some statements in \cite{NM} are not correct (see Remark \ref{R:C1}).

\section{The Briot-Bouquet Systems}\label{S:BB}

Since the work of Briot and Bouquet \cite{BB}, many authors have worked on the Briot-Bouquet systems (see \cite{IKSY} for an extensive bibliography). A \textit{Briot-Bouquet} system is of the following form
\begin{equation}\label{E:BB}
xy_i'=f_i(x,y_1,\cdots,y_n),\ \ \ \ \ \ (i=1,\cdots,n),
\end{equation}
where $f_i$'s are holomorphic and satisfy $f_i(0,0,\cdots,0)=0$.

We are interested in the existence of holomorphic solutions of (\ref{E:BB}) at $x=0$, for which the matrix
$$A=\left( \frac{\partial f_i}{\partial y_j}(0,0,\cdots,0) \right)_{1\le i,j\le n}$$
plays an important role. The one-dimensional theory is well-known since the original work of Briot and Bouquet (cf. \cite{H}). In higher dimensions, we have the following basic result (cf. \cite[\S 4, Proposition 1.1.1]{IKSY}).

\begin{prop}\label{P:A}
If none of the eigenvalues of $A$ is a positive integer, then \textup{(\ref{E:BB})} admits a unique holomorphic solution at $x=0$ satisfying $y_i(0)=0$.
\end{prop}

Surprisingly, it is very difficult to find an explicit treatment of the case when some of the eigenvalues of $A$ are positive integers in the existing literature. Therefore we give a detailed study of this case in this section. For simplicity, we only write out the details in dimension two, although it will be clear that similar results hold in higher dimensions.

Assume $n=2$. Let $q,s$ be the eigenvalues of $A$. First suppose that only $q$ is a positive integer. Then after a suitable linear change of variables (\ref{E:BB}) can be written as
\begin{equation}\label{E:1}
\left\{ \begin{aligned}
xu' & = px+qu+f(x,u,v),\\
xv' & = rx+sv+g(x,u,v),
\end{aligned} \right.
\end{equation}
where $f(x,u,v)=\sum_{|\alpha|\ge 2}a_\alpha x^{\alpha_1}u^{\alpha_2}v^{\alpha_3}$ and $g(x,u,v)=\sum_{|\beta|\ge 2}b_\beta x^{\beta_1}u^{\beta_2}v^{\beta_3}$.

Suppose that $u(x)=\sum_{k\ge 1}c_k x^k$, $v(x)=\sum_{k\ge 1}d_k x^k$ is a formal solution of (\ref{E:1}). Then it is easy to see that we have $(1-q)c_1=p$. Thus if $q=1$ and $p\neq 0$ then (\ref{E:1}) does not have any holomorphic solutions. On the other hand, if $q=1$ and $p=0$, then the choice of $c_1$ is arbitrary and other $c_k$ and $d_k$ can then be determined recursively. Moreover, as in the proof of Proposition \ref{P:A} in \cite{IKSY}, one can show that such formal solutions actually converge. Thus (\ref{E:1}) has infinitely many holomorphic solutions in this case.

If $q>1$, then we make the following change of variables:
$$u=x(\tilde{u}+\frac{p}{1-q}),\ \ \ v=x(\tilde{v}+\frac{r}{1-s}).$$
One readily checks that (\ref{E:1}) then takes the form
\begin{equation}\label{E:2}
\left\{ \begin{aligned}
x\tilde{u}' & = \tilde{p}x+(q-1)\tilde{u}+\tilde{f}(x,\tilde{u},\tilde{v}),\\
x\tilde{v}' & = \tilde{r}x+(s-1)\tilde{v}+\tilde{g}(x,\tilde{u},\tilde{v}),
\end{aligned} \right.
\end{equation}
where
$$\begin{aligned}
\tilde{p} & = a_{200}+\frac{p}{1-q}a_{110}+\frac{r}{1-s}a_{101}+\frac{pr}{(1-q)(1-s)}a_{011},\\
\tilde{r} & = b_{200}+\frac{p}{1-q}b_{110}+\frac{r}{1-s}b_{101}+\frac{pr}{(1-q)(1-s)}b_{011},
\end{aligned}$$
with $\tilde{f}(x,\tilde{u},\tilde{v})=\sum_{|\alpha|\ge 2}\tilde{a}_\alpha x^{\alpha_1}\tilde{u}^{\alpha_2}\tilde{v}^{\alpha_3}$ and $\tilde{g}(x,\tilde{u},\tilde{v})=\sum_{|\beta|\ge 2}\tilde{b}_\beta x^{\beta_1}\tilde{u}^{\beta_2}\tilde{v}^{\beta_3}$.

After $q-1$ such steps, (\ref{E:1}) then takes the form
\begin{equation}\label{E:3}
\left\{ \begin{aligned}
x\bar{u}' & = \bar{p}x+\bar{u}+\bar{f}(x,\bar{u},\bar{v}),\\
x\bar{v}' & = \bar{r}x+(s-q+1)\bar{v}+\bar{g}(x,\bar{u},\bar{v}).
\end{aligned} \right.
\end{equation}
Thus by the above discussion, (\ref{E:3}), hence (\ref{E:1}), has no holomorphic solutions if $\bar{p}\neq 0$ and infinitely many solutions if $\bar{p}=0$.

In summary, we have proved the following

\begin{prop}\label{P:1}
Assume that $A$ has one positive integer eigenvalue $q$. Let $\bar{p}$ be as in \textup{(\ref{E:3})}. If $\bar{p}\neq 0$ then \textup{(\ref{E:BB})} admits no holomorphic solutions at $x=0$. If $\bar{p}=0$ then \textup{(\ref{E:BB})} admits infinitely many holomorphic solutions at $x=0$ satisfying $y_i(0)=0$.
\end{prop}

Suppose now that both $q$ and $s$ are positive integers, with $q\le s$. First suppose that $A$ is diagonalizable. Then after a suitable linear change of variables (\ref{E:BB}) can again be written as in (\ref{E:1}). And (\ref{E:1}) takes the form (\ref{E:3}) after $q-1$ steps.

If $\bar{p}\neq 0$ then (\ref{E:3}) has no holomorphic solutions. Suppose $s=q$, then (\ref{E:3}) has no holomorphic solutions if $\bar{r}\neq 0$ and infinitely many holomorphic solutions if $\bar{p}=0$ and $\bar{r}=0$. Suppose $s>q$ and $\bar{p}=0$, then after another $s-q$ steps (\ref{E:1}) takes the form
\begin{equation}\label{E:4}
\left\{ \begin{aligned}
x\hat{u}' & = \hat{p}x+(q-s+1)\hat{u}+\hat{f}(x,\hat{u},\hat{v}),\\
x\hat{v}' & = \hat{r}x+\hat{v}+\hat{g}(x,\hat{u},\hat{v}).
\end{aligned} \right.
\end{equation}
Thus (\ref{E:4}) has no holomorphic solutions if $\hat{r}\neq 0$ and infinitely many holomorphic solutions if $\hat{r}=0$.

In summary, we have proved the following

\begin{prop}\label{P:2}
Assume that $A$ has two positive integer eigenvalues $q$ and $s$, with $q<s$. Let $\bar{p}$ be as in \textup{(\ref{E:3})} and $\hat{r}$ as in \textup{(\ref{E:4})}. If $\bar{p}\neq 0$ or $\bar{p}=0$ and $\hat{r}\neq 0$ then \textup{(\ref{E:BB})} admits no holomorphic solutions at $x=0$. If $\bar{p}=0$ and $\hat{r}=0$ then \textup{(\ref{E:BB})} admits infinitely many holomorphic solutions at $x=0$ satisfying $y_i(0)=0$.
\end{prop}

\begin{prop}\label{P:3}
Assume that $A$ has two equal positive integer eigenvalues $q$ and $A$ is diagonalizable. Let $\bar{p}$ and $\bar{r}$ be as in \textup{(\ref{E:3})}. If $\bar{p}\neq 0$ or $\bar{r}\neq 0$ then \textup{(\ref{E:BB})} admits no holomorphic solutions at $x=0$. If $\bar{p}=0$ and $\bar{r}=0$ then \textup{(\ref{E:BB})} admits infinitely many holomorphic solutions at $x=0$ satisfying $y_i(0)=0$.
\end{prop}

Finally suppose that $q=s$ and $A$ is not diagonalizable. Then after a suitable linear change of variables (\ref{E:BB}) can be written as
\begin{equation}\label{E:5}
\left\{ \begin{aligned}
xu' & = px+qu+\epsilon v+f(x,u,v),\\
xv' & = rx+qv+g(x,u,v).
\end{aligned} \right.
\end{equation}

Suppose first $q=1$. If $r\neq 0$ then (\ref{E:5}) has no holomorphic solutions. If $r=0$ then $d_1=-p/\epsilon$ and the choice of $c_1$ is arbitrary. Then (\ref{E:5}) has infinitely many holomorphic solutions in this case.

Suppose now $q>1$. We make the following type of change of variables:
$$u=x(\tilde{u}+\frac{p+\frac{\epsilon r}{1-q}}{1-q}),\ \ \ v=x(\tilde{v}+\frac{r}{1-q}).$$
One readily checks that after $q-1$ steps (\ref{E:5}) then takes the form
\begin{equation}\label{E:6}
\left\{ \begin{aligned}
x\bar{u}' & = \bar{p}x+\bar{u}+\epsilon \bar{v}+\bar{f}(x,\bar{u},\bar{v}),\\
x\bar{v}' & = \bar{r}x+\bar{v}+\bar{g}(x,\bar{u},\bar{v}).
\end{aligned} \right.
\end{equation}
Thus (\ref{E:6}) has no holomorphic solutions if $\bar{r}\neq 0$ and infinitely many holomorphic solutions if $\bar{r}=0$.

In summary, we have proved the following

\begin{prop}\label{P:4}
Assume that $A$ has two equal positive integer eigenvalues $q$ and $A$ is not diagonalizable. Let $\bar{r}$ be as in \textup{(\ref{E:6})}. If $\bar{r}\neq 0$ then \textup{(\ref{E:BB})} admits no holomorphic solutions at $x=0$. If $\bar{r}=0$ then \textup{(\ref{E:BB})} admits infinitely many holomorphic solutions at $x=0$ satisfying $y_i(0)=0$.
\end{prop}

\section{The Center Families}\label{S:CF}

We consider the three-dimensional complex system
\begin{equation}\label{E:C1}
\left\{ \begin{aligned}
\dot{x} & = F(x,y,z),\\
\dot{y} & = G(x,y,z),\\
\dot{z} & = H(x,y,z),
\end{aligned} \right.
\end{equation}
where $F$,$G$,$H$ are holomorphic and satisfy $F(0,0,0)=G(0,0,0)=H(0,0,0)=0$.

Denote by $\Lambda (x,y,z)^T$ the linear part of the above system. In \cite{N}, Needham proved the following (cf. \cite[Theorem 2]{N})

\begin{thm}\label{T:N}
Assume that $\Lambda$ has only one purely imaginary eigenvalue $i\omega$, with $x$-axis the eigenspace. Then \textup{(\ref{E:C1})} has a unique $x$-invariant holomorphic center manifold at the origin which contains an isochronous center family of period $2\pi/|\omega|$.
\end{thm}

Here ``$x$-invariant" means ``tangent to the $x$-axis". In what follows, we extend this result to the case when $\Lambda$ has more than one purely imaginary eigenvalues, based on our study of the Briot-Bouquet systems in the previous section.

First suppose that $\Lambda$ has two purely imaginary eigenvalues and $\Lambda$ is diagonalizable. By scaling, we assume that one of them is equal to $i$ and the other equal to $i\mu$ with $\mu\in \rv$ and $|\mu|\ge 1$. Then after a suitable linear change of variables, (\ref{E:C1}) can be written as
\begin{equation}\label{E:C2}
\left\{ \begin{aligned}
\dot{x} & = ix + f(x,y,z),\\
\dot{y} & = i\mu y + g(x,y,z),\\
\dot{z} & = \lambda z + h(x,y,z),
\end{aligned} \right.
\end{equation}
where $\textup{Re}\lambda\neq 0$, $f(x,y,z)=\sum_{|\alpha|\ge 2}a_\alpha x^{\alpha_1}y^{\alpha_2}z^{\alpha_3}$, $g(x,y,z)=\sum_{|\beta|\ge 2}b_\beta x^{\beta_1}y^{\beta_2}z^{\beta_3}$ and $h(x,y,z)=\sum_{|\gamma|\ge 2}c_\gamma x^{\gamma_1}y^{\gamma_2}z^{\gamma_3}$.

To find a center manifold of the form $(x(y),y,z(y))$, we consider the following initial value problem (cf. \cite{N},\cite{NM}).
\begin{equation}\label{E:C3}
\left\{ \begin{aligned}
(i\mu y+g(x(y),y,z(y))) x'(y) & = ix(y)+f(x(y),y,z(y)),\ \ \ x(0)=0,\\
(i\mu y+g(x(y),y,z(y))) z'(y) & = \lambda z(y)+h(x(y),y,z(y)),\ \ \ z(0)=0.
\end{aligned} \right.
\end{equation}

Set $x(y)=yu(y)$ and $z(y)=yv(y)$. Then we have
\begin{equation}\label{E:C4}
\left\{ \begin{aligned}
yu' & = (\frac{1}{\mu}-1)u+(\tilde{f}(u,y,v)-\frac{1}{\mu}u\tilde{g}(u,y,v))+\cdots,\\
yv' & = (-i\frac{\lambda}{\mu}-1)v+(\tilde{h}(u,y,v)+i\frac{\lambda}{\mu}v\tilde{g}(u,y,v))+\cdots,
\end{aligned} \right.
\end{equation}
where $\tilde{f}(u,y,v)=f(yu,y,yv)/(i\mu y)$, $\tilde{g}(u,y,v)=g(yu,y,yv)/(i\mu y)$ and $\tilde{h}(u,y,v)=h(yu,y,yv)/(i\mu y)$.

Since $-i\lambda/\mu-1\neq 0$, we have $v(0)=0$, i.e. $z'(0)=0$. If $\mu\neq 1$ then $1/\mu-1\neq 0$. Thus we have $u(0)=0$, i.e. $x'(0)=0$. By Proposition \ref{P:A}, (\ref{E:C4}) has a unique holomorphic solution at $y=0$. If $\mu=1$ then the choice of $u(0)$, i.e. $x'(0)$, is arbitrary. For each choice of $u(0)$, (\ref{E:C4}) has a unique holomorphic solution at $y=0$ by Proposition \ref{P:A}.

To find a center manifold of the form $(x,y(x),z(x))$, we consider the following initial value problem.
\begin{equation}\label{E:C5}
\left\{ \begin{aligned}
(ix+f(x,y(x),z(x))) y'(x) & = i\mu y(x)+g(x,y(x),z(x)),\ \ \ y(0)=0,\\
(ix+f(x,y(x),z(x))) z'(x) & = \lambda z(x)+h(x,y(x),z(x)),\ \ \ \ z(0)=0.
\end{aligned} \right.
\end{equation}

Set $y(x)=xu(x)$ and $z(x)=xv(x)$. Then we have
\begin{equation}\label{E:C6}
\left\{ \begin{aligned}
xu' & = (\mu-1)u+(\tilde{g}(x,u,v)-\mu u\tilde{f}(x,u,v))+\cdots,\\
xv' & = (-i\lambda-1)v+(\tilde{h}(x,u,v)+i\lambda v\tilde{f}(x,u,v))+\cdots,
\end{aligned} \right.
\end{equation}
where $\tilde{f}(x,u,v)=f(x,xu,xv)/(ix)$, $\tilde{g}(x,u,v)=g(x,xu,xv)/(ix)$ and $\tilde{h}(x,u,v)=h(x,xu,xv)/(ix)$.

Since $-i\lambda-1\neq 0$, we have $v(0)=0$, i.e. $z'(0)=0$. If $\mu\neq 1$ then $\mu-1\neq 0$. Thus we have $u(0)=0$, i.e. $y'(0)=0$. Thus if $\mu\not\in \nv$ then (\ref{E:C6}) has a unique holomorphic solution at $x=0$ by Proposition \ref{P:A}. If $\mu=1$ then the choice of $u(0)$, i.e. $y'(0)$, is arbitrary. For each choice of $u(0)$, (\ref{E:C6}) has a unique holomorphic solution at $y=0$ by Proposition \ref{P:A}. If $\mu\in\nv\backslash \{1\}$, then we rewrite (\ref{E:C6}) as
\begin{equation}\label{E:C7}
\left\{ \begin{aligned}
xu' & = (\mu-1)u-ib_{200}x+O(2),\\
xv' & = (-i\lambda-1)v-ic_{200}x+O(2).
\end{aligned} \right.
\end{equation}

If $\mu=2$ then, by Proposition \ref{P:1}, (\ref{E:C7}) has no holomorphic solutions if $b_{200}\neq 0$ and infinitely many holomorphic solutions if $b_{200}=0$. (In the latter case we always have $u(0)=0$, i.e. $y'(0)=0$, but the choice of $u'(0)$, i.e. $y''(0)$, is arbitrary.) If $\mu>2$ then we make $\mu-2$ steps of change of variables as in the previous section, after which (\ref{E:C7}) becomes
\begin{equation}\label{E:C8}
\left\{ \begin{aligned}
x\bar{u}' & = \bar{u}+\bar{p}x+O(2),\\
x\bar{v}' & = (-i\lambda-\mu+1)\bar{v}+\bar{r}x+O(2),
\end{aligned} \right.
\end{equation}
where $\bar{p}$ (resp. $\bar{r}$) is determined by the coefficients of $f$ and $g$ (resp. $h$). We can then again apply Proposition \ref{P:1}.

Note that each holomorphic solution of the Briot-Bouquet system gives a center manifold of the original complex dynamical system. On each such center manifold there exists an isochronous center family (due to a holomorphic linearization of the system on the center manifold, cf. \cite{A}).

In summary, we have proved the following

\begin{thm}\label{T:C1}
Assume that $\Lambda$ has two distinct purely imaginary eigenvalues $i\omega_1$ and $i\omega_2$, with $x$-axis and $y$-axis the respective eigenspaces, and $|\omega_1|\le|\omega_2|$. Then

$1)$ System \textup{(\ref{E:C1})} has a $y$-invariant holomorphic center manifold at the origin which contains an isochronous center family of period $2\pi/|\omega_2|$.

$2)$ If $\omega_2/\omega_1\not\in\nv$ then \textup{(\ref{E:C1})} has a $x$-invariant holomorphic center manifold at the origin which contains an isochronous center family of period $2\pi/|\omega_1|$.

$3)$ If $\omega_2/\omega_1\in\nv$, let $\bar{p}$ be as in \textup{(\ref{E:C8})}.

$3.1)$ If $\bar{p}=0$ then \textup{(\ref{E:C1})} has infinitely many $x$-invariant holomorphic center manifolds at the origin, each of which contains an isochronous center family of period $2\pi/|\omega_1|$.

$3.2)$ If $\bar{p}\neq 0$ then \textup{(\ref{E:C1})} has no other holomorphic center manifolds at the origin.
\end{thm}

\begin{thm}\label{T:C2}
Assume that $\Lambda$ has two equal purely imaginary eigenvalues $i\omega$, with $(x,y)$-plane the eigenspace, and $\Lambda$ is diagonalizable. Then \textup{(\ref{E:C1})} has a unique $(x,y)$-invariant holomorphic center manifold at the origin which contains an isochronous center family of period $2\pi/|\omega|$.
\end{thm}

\begin{rmk}\label{R:C1}
The statement of Theorem 6.1(c) in \cite{NM} is incorrect. First of all the two different subcases there are not determined by $b_{\sigma 0}$ but rather by $\bar{p}$ as in the above theorem. Second, when $b_{\sigma 0}=0$ (which should be $\bar{p}=0$), $(z,w)=(0,0)$ (as $(x,y)=(0,0)$ in our case) is not an isochronous center. Because of this, some other statements (such as those in section 7 of \cite{NM}) also need to be modified accordingly.
\end{rmk}

Now suppose that $\Lambda$ has two equal purely imaginary eigenvalues and $\Lambda$ is not diagonalizable. After scaling and a suitable linear change of variables, (\ref{E:C1}) can be written as
\begin{equation}\label{E:C9}
\left\{ \begin{aligned}
\dot{x} & = ix + y + f(x,y,z),\\
\dot{y} & = iy + g(x,y,z),\\
\dot{z} & = \lambda z + h(x,y,z).
\end{aligned} \right.
\end{equation}

To find a center manifold of the form $(x(y),y,z(y))$, we consider the following initial value problem.
\begin{equation}\label{E:C10}
\left\{ \begin{aligned}
(iy+g(x(y),y,z(y))) x'(y) & = ix(y)+y+f(x(y),y,z(y)),\ \ x(0)=0,\\
(iy+g(x(y),y,z(y))) z'(y) & = \lambda z(y)+h(x(y),y,z(y)),\ \ \ \ \ \ \ z(0)=0.
\end{aligned} \right.
\end{equation}

Set $x(y)=yu(y)$ and $z(y)=yv(y)$. Then we have
\begin{equation}\label{E:C11}
\left\{ \begin{aligned}
yu' & = -i+O(1),\\
yv' & = (-i\lambda-1)v+\cdots.
\end{aligned} \right.
\end{equation}

Due to the presence of $-i$, (\ref{E:C11}) has no solutions.

To find a center manifold of the form $(x,y(x),z(x))$, we consider the following initial value problem.
\begin{equation}\label{E:C12}
\left\{ \begin{aligned}
(ix+y(x)+f(x,y(x),z(x))) y'(x) & = iy(x)+g(x,y(x),z(x)),\ \ y(0)=0,\\
(ix+y(x)+f(x,y(x),z(x))) z'(x) & = \lambda z(x)+h(x,y(x),z(x)),\ z(0)=0.
\end{aligned} \right.
\end{equation}

Set $y(x)=xu(x)$ and $z(x)=xv(x)$. Then we have
\begin{equation}\label{E:C13}
\left\{ \begin{aligned}
xu' & = 0u+0v+\cdots,\\
xv' & = iu+(-i\lambda-1)v+\cdots.
\end{aligned} \right.
\end{equation}

Since there is no center manifold of the form $(x(y),y,z(y))$, we must have $y'(0)=0$, i.e. $u(0)=0$, which implies $v(0)=0$, i.e. $z'(0)=0$. Then by Proposition \ref{P:A}, (\ref{E:C13}) has a unique holomorphic solution at $x=0$.

In summary, we have proved the following

\begin{thm}\label{T:C3}
Assume that $\Lambda$ has two equal purely imaginary eigenvalues $i\omega$, with $(x,y)$-plane the eigenspace, and $\Lambda$ is not diagonalizable. Then \textup{(\ref{E:C1})} has a unique $x$-invariant holomorphic center manifold at the origin which contains an isochronous center family of period $2\pi/|\omega|$.
\end{thm}

\begin{rmk}\label{R:C2}
In \cite{NM}, Needham and McAllister studied the two dimensional case but assuming that $\Lambda$ is diagonalizable. If we assume that $\Lambda$ is not diagonalizable, then a similar result as the above theorem holds for the two dimensional case.
\end{rmk}

Next suppose that $\Lambda$ has three purely imaginary eigenvalues and $\Lambda$ is diagonalizable. By scaling, we assume that one of them is equal to $i$ and the other two equal to $i\mu$ and $i\nu$ with $\mu,\nu\in \rv$ and $|\nu|\ge|\mu|\ge 1$. Then after a suitable linear change of variables, (\ref{E:C1}) can be written as
\begin{equation}\label{E:C14}
\left\{ \begin{aligned}
\dot{x} & = ix + f(x,y,z),\\
\dot{y} & = i\mu y + g(x,y,z),\\
\dot{z} & = i\nu z + h(x,y,z).
\end{aligned} \right.
\end{equation}

To find a center manifold of the form $(x(z),y(z),z)$, we consider the following initial value problem.
\begin{equation}\label{E:C15}
\left\{ \begin{aligned}
(i\nu z+h(x(z),y(z),z)) x'(z) & = ix(z)+f(x(z),y(z),z),\ \ \ \ x(0)=0,\\
(i\nu z+h(x(z),y(z),z)) y'(z) & = i\mu y(z)+g(x(z),y(z),z),\ \ \ y(0)=0.
\end{aligned} \right.
\end{equation}

Set $x(z)=zu(z)$ and $y(z)=zv(z)$. Then we have
\begin{equation}\label{E:C16}
\left\{ \begin{aligned}
zu' & = (\frac{1}{\nu}-1)u+(\tilde{f}(u,v,z)-\frac{1}{\nu}u\tilde{h}(u,v,z))+\cdots,\\
zv' & = (\frac{\mu}{\nu}-1)v+(\tilde{g}(u,v,z)-\frac{\mu}{\nu}v\tilde{h}(u,v,z))+\cdots,
\end{aligned} \right.
\end{equation}
where $\tilde{f}(u,v,z)=f(zu,zv,z)/(i\nu z)$, $\tilde{g}(u,v,z)=g(zu,zv,z)/(i\nu z)$ and $\tilde{h}(u,v,z)=h(zu,zv,z)/(i\nu z)$.

The choice of $u(0)$, i.e. $x'(0)$, is arbitrary if and only if $\nu=1$ (and $\mu=\pm 1$). The choice of $v(0)$, i.e. $y'(0)$, is arbitrary if and only if $\mu=\nu$. For each possible choice of $(u(0),v(0))$, (\ref{E:C16}) has a unique holomorphic solution at $z=0$ by Proposition \ref{P:A}.

To find a center manifold of the form $(x(y),y,z(y))$, we consider the following initial value problem.
\begin{equation}\label{E:C17}
\left\{ \begin{aligned}
(i\mu y+g(x(y),y,z(y))) x'(y) & = ix(y)+f(x(y),y,z(y)),\ \ \ \ x(0)=0,\\
(i\mu y+g(x(y),y,z(y))) z'(y) & = i\nu z(y)+h(x(y),y,z(y)),\ \ \ z(0)=0.
\end{aligned} \right.
\end{equation}

Set $x(y)=yu(y)$ and $z(y)=yv(y)$. Then we have
\begin{equation}\label{E:C18}
\left\{ \begin{aligned}
yu' & = (\frac{1}{\mu}-1)u+(\tilde{f}(u,y,v)-\frac{1}{\mu}u\tilde{g}(u,y,v))+\cdots,\\
yv' & = (\frac{\nu}{\mu}-1)v+(\tilde{h}(u,y,v)-\frac{\nu}{\mu}v\tilde{g}(u,y,v))+\cdots,
\end{aligned} \right.
\end{equation}
where $\tilde{f}(u,y,v)=f(yu,y,yv)/(i\mu y)$, $\tilde{g}(u,y,v)=g(yu,y,yv)/(i\mu y)$ and $\tilde{h}(u,y,v)=h(yu,y,yv)/(i\mu y)$.

The choice of $u(0)$, i.e. $x'(0)$, is arbitrary if and only if $\mu=1$. The choice of $v(0)$, i.e. $z'(0)$, is arbitrary if and only if $\mu=\nu$. If $\nu/\mu\not\in\nv\backslash \{1\}$ then for each possible choice of $(u(0),v(0))$, (\ref{E:C18}) has a unique holomorphic solution at $y=0$ by Proposition \ref{P:A}. If $\nu/\mu\in\nv\backslash \{1\}$, then we rewrite (\ref{E:C18}) as
\begin{equation}\label{E:C19}
\left\{ \begin{aligned}
yu' & = (\frac{1}{\mu}-1)u-ia_{020}y+O(2),\\
yv' & = (\frac{\nu}{\mu}-1)v-ic_{020}y+O(2).
\end{aligned} \right.
\end{equation}

If $\nu/\mu=2$ then, by Proposition \ref{P:1}, (\ref{E:C19}) has no holomorphic solutions if $c_{020}\neq 0$ and infinitely many holomorphic solutions if $c_{020}=0$. (In the latter case we always have $v(0)=0$, i.e. $z'(0)=0$, but the choice of $v'(0)$, i.e. $z''(0)$, is arbitrary.) If $\nu/\mu>2$ then we make $\nu/\mu-2$ steps of change of variables as in the previous section, after which (\ref{E:C19}) becomes
\begin{equation}\label{E:C20}
\left\{ \begin{aligned}
y\bar{u}' & = (\frac{1-\nu}{\mu}+1)\bar{u}+\bar{p}y+O(2),\\
y\bar{v}' & = \bar{v}+\bar{r}y+O(2),
\end{aligned} \right.
\end{equation}
where $\bar{p}$ (resp. $\bar{r}$) is determined by the coefficients of $f$ (resp. $h$) and $g$. We can then again apply Proposition \ref{P:1}.

To find a center manifold of the form $(x,y(x),z(x))$, we consider the following initial value problem.
\begin{equation}\label{E:C21}
\left\{ \begin{aligned}
(ix+f(x,y(x),z(x))) y'(x) & = i\mu y(x)+g(x,y(x),z(x)),\ \ \ y(0)=0,\\
(ix+f(x,y(x),z(x))) z'(x) & = i\nu z(x)+h(x,y(x),z(x)),\ \ \ z(0)=0.
\end{aligned} \right.
\end{equation}

Set $y(x)=xu(x)$ and $z(x)=xv(x)$. Then we have
\begin{equation}\label{E:C22}
\left\{ \begin{aligned}
xu' & = (\mu-1)u+(\tilde{g}(x,u,v)-\mu u\tilde{f}(x,u,v))+\cdots,\\
xv' & = (\nu-1)v+(\tilde{h}(x,u,v)-\nu v\tilde{f}(x,u,v))+\cdots,
\end{aligned} \right.
\end{equation}
where $\tilde{f}(x,u,v)=f(x,xu,xv)/(ix)$, $\tilde{g}(x,u,v)=g(x,xu,xv)/(ix)$ and $\tilde{h}(x,u,v)=h(x,xu,xv)/(ix)$.

The choice of $u(0)$, i.e. $y'(0)$, is arbitrary if and only if $\mu=1$. The choice of $v(0)$, i.e. $z'(0)$, is arbitrary if and only if $\nu=1$ (and $\mu=\pm 1$). If neither $\mu$ nor $\nu$ belongs to $\nv\backslash \{1\}$ then for each possible choice of $(u(0),v(0))$, (\ref{E:C22}) has a unique holomorphic solution at $x=0$ by Proposition \ref{P:A}. If $\mu\in \nv\backslash \{1\}$ but $\nu\not\in \nv\backslash \{1\}$ then after making $\mu-2$ steps of change of variables as in the previous section, (\ref{E:C22}) becomes
\begin{equation}\label{E:C23}
\left\{ \begin{aligned}
x\bar{u}' & = \bar{u}+\bar{p}x+O(2),\\
x\bar{v}' & = (\nu-\mu+1)\bar{v}+\bar{r}x+O(2),
\end{aligned} \right.
\end{equation}
where $\bar{p}$ (resp. $\bar{r}$) is determined by the coefficients of $f$ and $g$ (resp. $h$). By Proposition \ref{P:1}, (\ref{E:C23}) has no holomorphic solutions at $x=0$ if $\bar{p}\neq 0$ and infinitely many holomorphic solutions at $x=0$ if $\bar{p}=0$. (In the latter case the choice of $\bar{u}'(0)$, i.e. $y^{(\mu)}(0)$, is arbitrary.)

If $\nu\in \nv\backslash \{1\}$ but $\mu\not\in \nv\backslash \{1\}$ then after making $\nu-2$ steps of change of variables as in the previous section, (\ref{E:C22}) becomes
\begin{equation}\label{E:C24}
\left\{ \begin{aligned}
x\hat{u}' & = (\mu-\nu+1)\hat{u}+\hat{p}x+O(2),\\
x\hat{v}' & = \hat{v}+\hat{r}x+O(2),
\end{aligned} \right.
\end{equation}
where $\hat{p}$ (resp. $\hat{r}$) is determined by the coefficients of $f$ and $g$ (resp. $h$). By Proposition \ref{P:1}, (\ref{E:C24}) has no holomorphic solutions at $x=0$ if $\hat{r}\neq 0$ and infinitely many holomorphic solutions at $x=0$ if $\hat{r}=0$. (In the latter case the choice of $\hat{v}'(0)$, i.e. $z^{(\nu)}(0)$, is arbitrary.)

Assume now that both $\mu$ and $\nu$ belong to $\nv\backslash \{1\}$. If $\nu=\mu$ then, by Proposition \ref{P:3}, (\ref{E:C22}) has no holomorphic solutions at $x=0$ if $\bar{p}\neq 0$ or $\bar{r}\neq 0$ and infinitely many holomorphic solutions at $x=0$ if $\bar{p}=0$ and $\bar{r}=0$. (In the latter case the choices of $\bar{u}'(0)$ and $\bar{v}'(0)$, i.e. $y^{(\mu)}(0)$ and $z^{(\mu)}(0)$, are arbitrary.)

If $\nu>\mu$ then, by Proposition \ref{P:2}, (\ref{E:C22}) has no holomorphic solutions at $x=0$ if $\bar{p}\neq 0$ or $\bar{p}=0$ and $\hat{r}\neq 0$ and infinitely many holomorphic solutions at $x=0$ if $\bar{p}=0$ and $\hat{r}=0$. (In the latter case the choices of $\bar{u}'(0)$ and $\hat{v}'(0)$, i.e. $y^{(\mu)}(0)$ and $z^{(\nu)}(0)$, are arbitrary.)

In summary, we have proved the following

\begin{thm}\label{T:C4}
Assume that $\Lambda$ has three distinct purely imaginary eigenvalues $i\omega_1$, $i\omega_2$ and $i\omega_3$, with $x$-axis, $y$-axis and $z$-axis the respective eigenspaces, and $|\omega_1|\le|\omega_2|\le|\omega_3|$. Then

$1)$ System \textup{(\ref{E:C1})} has a $z$-invariant holomorphic center manifold at the origin which contains an isochronous center family of period $2\pi/|\omega_3|$.

$2)$ If $\omega_3/\omega_2\not\in\nv$ then \textup{(\ref{E:C1})} has a $y$-invariant holomorphic center manifold at the origin which contains an isochronous center family of period $2\pi/|\omega_2|$.

$3)$ If $\omega_3/\omega_2\in\nv$, let $\bar{p}$ be as in \textup{(\ref{E:C20})}.

$3.1)$ If $\bar{p}=0$ then \textup{(\ref{E:C1})} has infinitely many $y$-invariant holomorphic center manifolds at the origin, each of which contains an isochronous center family of period $2\pi/|\omega_2|$.

$3.2)$ If $\bar{p}\neq 0$ then \textup{(\ref{E:C1})} has no $y$-invariant holomorphic center manifolds at the origin. 

$4)$ If $\omega_2/\omega_1\not\in\nv$ and $\omega_3/\omega_1\not\in\nv$ then \textup{(\ref{E:C1})} has a $x$-invariant holomorphic center manifold at the origin which contains an isochronous center family of period $2\pi/|\omega_1|$.

$5)$ If $\omega_2/\omega_1\in\nv$ but $\omega_3/\omega_1\not\in\nv$, let $\bar{p}$ be as in \textup{(\ref{E:C23})}.

$5.1)$ If $\bar{p}=0$ then \textup{(\ref{E:C1})} has infinitely many $x$-invariant holomorphic center manifolds at the origin, each of which contains an isochronous center family of period $2\pi/|\omega_1|$.

$5.2)$ If $\bar{p}\neq 0$ then \textup{(\ref{E:C1})} has no $x$-invariant holomorphic center manifolds at the origin.

$6)$ If $\omega_2/\omega_1\not\in\nv$ but $\omega_3/\omega_1\in\nv$, let $\hat{r}$ be as in \textup{(\ref{E:C24})}.

$6.1)$ If $\hat{r}=0$ then \textup{(\ref{E:C1})} has infinitely many $x$-invariant holomorphic center manifolds at the origin, each of which contains an isochronous center family of period $2\pi/|\omega_1|$.

$6.2)$ If $\hat{r}\neq 0$ then \textup{(\ref{E:C1})} has no $x$-invariant holomorphic center manifolds at the origin.

$7)$ If $\omega_2/\omega_1\in\nv$ and $\omega_3/\omega_1\in\nv$, let $\bar{p}$ be as in \textup{(\ref{E:C23})} and $\hat{r}$ as in \textup{(\ref{E:C24})}.

$7.1)$ If $\bar{p}=0$ and $\hat{r}=0$ then \textup{(\ref{E:C1})} has infinitely many $x$-invariant holomorphic center manifolds at the origin, each of which contains an isochronous center family of period $2\pi/|\omega_1|$.

$7.2)$ If $\bar{p}\neq 0$ or $\bar{p}=0$ and $\hat{r}\neq 0$ then \textup{(\ref{E:C1})} has no $x$-invariant holomorphic center manifolds at the origin.
\end{thm}

\begin{thm}\label{T:C5}
Assume that $\Lambda$ has two equal purely imaginary eigenvalues $i\omega_1$ and another purely imaginary eigenvalue $i\omega_2$, with $(x,y)$-plane and $z$-axis the respective eigenspaces, $\omega_1\neq \omega_2$ and $|\omega_1|\le|\omega_2|$ and $\Lambda$ is diagonalizable. Then

$1)$ System \textup{(\ref{E:C1})} has a $z$-invariant holomorphic center manifold at the origin which contains an isochronous center family of period $2\pi/|\omega_2|$.

$2)$ Let $\bar{r}$ be as in \textup{(\ref{E:C20})} and $\hat{r}$ as in \textup{(\ref{E:C24})}.

$2.1)$ If $\omega_2/\omega_1\not\in\nv$ or $\omega_2/\omega_1\in\nv$ and $\bar{r}=0=\hat{r}$ then \textup{(\ref{E:C1})} has infinitely many $(x,y)$-invariant holomorphic center manifolds at the origin, each of which contains an isochronous center family of period $2\pi/|\omega_1|$.

$2.2)$ If $\omega_2/\omega_1\in\nv$ and $\bar{r}=0\neq\hat{r}$ then \textup{(\ref{E:C1})} has infinitely many $y$-invariant holomorphic center manifolds at the origin, each of which contains an isochronous center family of period $2\pi/|\omega_1|$.

$2.3)$ If $\omega_2/\omega_1\in\nv$ and $\bar{r}\neq 0=\hat{r}$ then \textup{(\ref{E:C1})} has infinitely many $x$-invariant holomorphic center manifolds at the origin, each of which contains an isochronous center family of period $2\pi/|\omega_1|$.
\end{thm}

\begin{thm}\label{T:C6}
Assume that $\Lambda$ has one purely imaginary eigenvalue $i\omega_1$ and two other equal purely imaginary eigenvalues $i\omega_2$, with $x$-axis and $(y,z)$-plane the respective eigenspaces, $\omega_1\neq \omega_2$ and $|\omega_1|\le|\omega_2|$ and $\Lambda$ is diagonalizable. Then 

$1)$ System \textup{(\ref{E:C1})} has a $(y,z)$-invariant holomorphic center manifold at the origin which contains an isochronous center family of period $2\pi/|\omega_2|$.

$2)$ If $\omega_2/\omega_1\not\in\nv$ then \textup{(\ref{E:C1})} has a $x$-invariant holomorphic center manifold at the origin which contains an isochronous center family of period $2\pi/|\omega_1|$.

$3)$ If $\omega_2/\omega_1\in\nv$, let $\bar{p}$ and $\bar{r}$ be as in \textup{(\ref{E:C23})}.

$3.1)$ If $\bar{p}=0$ and $\bar{r}=0$ then \textup{(\ref{E:C1})} has infinitely many $x$-invariant holomorphic center manifolds at the origin, each of which contains an isochronous center family of period $2\pi/|\omega_1|$.

$3.2)$ If $\bar{p}\neq 0$ or $\bar{r}\neq 0$ then \textup{(\ref{E:C1})} has no $x$-invariant holomorphic center manifolds at the origin. 
\end{thm}

\begin{thm}[Poincar\'{e} Isochronous Center Theorem]\label{T:C7}
Assume that $\Lambda$ has three equal purely imaginary eigenvalues $i\omega$ and $\Lambda$ is diagonalizable. Then \textup{(\ref{E:C1})} has an isochronous center at the origin, with period $2\pi/|\omega|$.
\end{thm}

Finally suppose that $\Lambda$ has three purely imaginary eigenvalues and $\Lambda$ is not diagonalizable. If $\Lambda$ has a $2\times 2$ Jordan block then after scaling and a suitable linear change of variables, (\ref{E:C1}) can be written as in (\ref{E:C9}), with $\lambda$ replaced by $i\mu$. Again due to the presence of $-i$ in (\ref{E:C11}) (with $\lambda$ replaced by $i\mu$), there is no center manifolds of the form $(x(y),y,z(y))$. And if $\mu\not\in \nv\backslash \{1\}$ then, for each possible choice of $z'(0)$ ($y'(0)=0$), we have a unique holomorphic center manifold of the form $(x,y(x),z(x))$. 

If $\mu\in \nv\backslash \{1\}$ then after $\mu-2$ steps of change of variables similarly as in the previous section, we get
\begin{equation}\label{E:C25}
\left\{ \begin{aligned}
x\tilde{u}' & = (2-\mu)\tilde{u}+\tilde{p}x+O(2),\\
x\tilde{v}' & = i\tilde{u}+\tilde{v}+\tilde{r}x+O(2).
\end{aligned} \right.
\end{equation}
After a suitable linear change of variables, we then have
\begin{equation}\label{E:C26}
\left\{ \begin{aligned}
x\bar{u}' & = (2-\mu)\bar{u}+\bar{p}x+O(2),\\
x\bar{v}' & = \bar{v}+\bar{r}x+O(2).
\end{aligned} \right.
\end{equation}
By Proposition \ref{P:1}, (\ref{E:C26}) has no holomorphic solutions at $x=0$ if $\bar{r}\neq 0$ and infinitely many holomorphic solutions at $x=0$ if $\bar{r}=0$.

To find a center manifold of the form $(x(z),y(z),z)$, we consider the following initial value problem.
\begin{equation}\label{E:C27}
\left\{ \begin{aligned}
(i\mu z+h(x(z),y(z),z)) x'(z) & = ix(z)+y(z)+f(x(z),y(z),z),\ x(0)=0,\\
(i\mu z+h(x(z),y(z),z)) y'(z) & = iy(z)+g(x(z),y(z),z),\ \ \ \ \ \ \ \ \ \ \ y(0)=0.
\end{aligned} \right.
\end{equation}

Set $x(z)=zu(z)$ and $y(z)=zv(z)$. Then we have
\begin{equation}\label{E:C28}
\left\{ \begin{aligned}
zu' & = (\frac{1}{\mu}-1)u-\frac{i}{\mu}v+(\tilde{f}(u,v,z)-\frac{1}{\mu}(u-iv)\tilde{h}(u,v,z))+\cdots,\\
zv' & = (\frac{1}{\mu}-1)v+(\tilde{g}(u,v,z)-\frac{1}{\mu}v\tilde{h}(u,v,z))+\cdots,
\end{aligned} \right.
\end{equation}
where $\tilde{f}(u,v,z)=f(zu,zv,z)/(i\mu z)$, $\tilde{g}(u,v,z)=g(zu,zv,z)/(i\mu z)$ and $\tilde{h}(u,v,z)=h(zu,zv,z)/(i\mu z)$.

If $1/\mu\not\in\nv\backslash \{1\}$ then, for each possible choice of $u(0)$ ($v(0)=0$), (\ref{E:C28}) has a unique holomorphic solution at $z=0$ by Proposition \ref{P:A}. If $1/\mu\in\nv\backslash \{1\}$, then after $1/\mu-2$ steps of change of variables similarly as in the previous section, we get
\begin{equation}\label{E:C29}
\left\{ \begin{aligned}
z\bar{u}' & = \bar{u}-\frac{i}{\mu}\bar{v}+\bar{p}z+O(2),\\
z\bar{v}' & = \bar{v}+\bar{r}z+O(2).
\end{aligned} \right.
\end{equation}
By Proposition \ref{P:4}, (\ref{E:C29}) has no holomorphic solutions at $z=0$ if $\bar{r}\neq 0$ and infinitely many holomorphic solutions at $z=0$ if $\bar{r}=0$.

If $\Lambda$ has a $3\times 3$ Jordan block then after scaling and a suitable linear change of variables, (\ref{E:C1}) can be written as
\begin{equation}\label{E:C30}
\left\{ \begin{aligned}
\dot{x} & = ix + y + f(x,y,z),\\
\dot{y} & = iy + z + g(x,y,z),\\
\dot{z} & = iz + h(x,y,z).
\end{aligned} \right.
\end{equation}

Similar to the $2\times 2$ Jordan block case, there are no center manifolds of the form $(x(y),y,z(y))$ or $(x(z),y(z),z)$ (due to the presence of $-i$). To find a center manifold of the form $(x,y(x),z(x))$, we consider the following initial value problem.
\begin{equation}\label{E:C31}
\left\{ \begin{aligned}
(ix+y(x)+f(x,y(x),z(x))) y'(x) & = iy(x)+z(x)+g(x,y(x),z(x)),\ y(0)=0,\\
(ix+y(x)+f(x,y(x),z(x))) z'(x) & = iz(x)+h(x,y(x),z(x)),\ \ \ \ \ \ \ \ \ \ z(0)=0.
\end{aligned} \right.
\end{equation}

Set $y(x)=xu(x)$ and $z(x)=xv(x)$. Then we have
\begin{equation}\label{E:C32}
\left\{ \begin{aligned}
xu' & = 0u-iv+\cdots,\\
xv' & = 0u+0v+\cdots.
\end{aligned} \right.
\end{equation}
By Proposition \ref{P:A}, (\ref{E:C32}) has a unique holomorphic solution at $x=0$ (with $u(0)=0$, i.e. $y'(0)=0$, and $v(0)=0$, i.e. $z'(0)=0$).

In summary, we have proved the following

\begin{thm}\label{T:C8}
Assume that $\Lambda$ has two equal purely imaginary eigenvalues $i\omega_1$ and another purely imaginary eigenvalue $i\omega_2$, with $(x,y)$-plane and $z$-axis the respective eigenspaces, and $\Lambda$ is not diagonalizable. Then

$1)$ System \textup{(\ref{E:C1})} has no $y$-invariant holomorphic center manifolds at the origin.

$2)$ If $\omega_2/\omega_1\not\in\nv$ then \textup{(\ref{E:C1})} has a unique $x$-invariant holomorphic center manifold at the origin which contains an isochronous center family of period $2\pi/|\omega_1|$.

$3)$ If $\omega_1/\omega_2\not\in\nv$ then \textup{(\ref{E:C1})} has a unique $z$-invariant holomorphic center manifold at the origin which contains an isochronous center family of period $2\pi/|\omega_2|$.

$4)$ If $\omega_1=\omega_2$ then \textup{(\ref{E:C1})} has a unique $(x,z)$-invariant holomorphic center manifold at the origin which contains an isochronous center family of period $2\pi/|\omega_1|$.

$5)$ If $\omega_2/\omega_1\in\nv\backslash \{1\}$, let $\bar{r}$ be as in \textup{(\ref{E:C26})}.

$5.1)$ If $\bar{r}=0$ then \textup{(\ref{E:C1})} has infinitely many $x$-invariant holomorphic center manifolds at the origin, each of which contains an isochronous center family of period $2\pi/|\omega_1|$.

$5.2)$ If $\bar{r}\neq 0$ then \textup{(\ref{E:C1})} has no $x$-invariant holomorphic center manifolds at the origin.

$6)$ If $\omega_1/\omega_2\in\nv\backslash \{1\}$, let $\bar{r}$ be as in \textup{(\ref{E:C29})}.

$6.1)$ If $\bar{r}=0$ then \textup{(\ref{E:C1})} has infinitely many $z$-invariant holomorphic center manifolds at the origin, each of which contains an isochronous center family of period $2\pi/|\omega_2|$.

$6.2)$ If $\bar{r}\neq 0$ then \textup{(\ref{E:C1})} has no $z$-invariant holomorphic center manifolds at the origin. 
\end{thm}

\begin{thm}\label{T:C9}
Assume that $\Lambda$ has three equal purely imaginary eigenvalues $i\omega$ and $\Lambda$ has a $3\times 3$ Jordan block. Then 

$1)$ System \textup{(\ref{E:C1})} has a unique $x$-invariant holomorphic center manifold at the origin which contains an isochronous center family of period $2\pi/|\omega|$.

$2)$ System \textup{(\ref{E:C1})} has no $y$-invariant or $z$-invariant holomorphic center manifolds at the origin.
\end{thm}

\begin{rmk}
With the same method, one can readily generalize all the above results to higher dimensions.
\end{rmk}


\begin{thebibliography}{99}

\bibitem{A}
Arnold, V.I.,
\emph{Geometric Methods in the Theory of Ordinary Differential Equations}, 2nd. ed., Springer, Berlin, 1988.

\bibitem{BB}
Briot, C., Bouquet, J.C.,
\emph{Recherches sur les propri\'{e}t\'{e}s des fonctions d\'{e}finies par des \'{e}quations diff\'{e}rentielles}, J. \'{E}col. Imp. Poly. {\bf 21} (1856), 133-197.

\bibitem{H}
Hille, E,
\emph{Ordinary Differential Equations in the Complex Domain}, John Wiley \& Sons, Inc., New York, 1976.

\bibitem{IKSY}
Iwasaki, K., Kimura, H., Shimomura, S., Yoshida, M.,
\emph{From Gauss to Painlev\'{e}. A Mordern Theory of Special Functions}, Aspects of Mathematics, Vol. {\bf E16}, Vieweg, Braunschweig, 1991.

\bibitem{N}
Needham, D.J.,
\emph{A centre theorem for two-dimensional complex holomorphic systems and its generalization}, Proc. R. Soc. Lond. A {\bf 450} (1995), 225-232.

\bibitem{NM}
Needham, D.J., McAllister, S.,
\emph{Centre families in two-dimensional complex holomorphic dynamical systems}, Proc. R. Soc. Lond. A {\bf 454} (1998), 2267-2278.

\bibitem{Z1}
Zhang, G.Y.,
\emph{Fixed point indices and invariant periodic sets of holomorphic systems}, Proc. Amer. Math. Soc. {\bf 135} (2007), 767–776.

\bibitem{Z2}
Zhang, G.Y.,
\emph{From an iteration formula to Poincar\'{e}'s isochronous center theorem for holomorphic vector fields}, Proc. Amer. Math. Soc. {\bf 135} (2007), 2887–2891.

\end{thebibliography}
\end{document}